\documentclass[12pt]{article}
\usepackage[pdftex]{color, graphicx}
\usepackage{setspace}
\usepackage[utf8]{inputenc}
\usepackage{amsmath}
\usepackage[pdftex,bookmarks,colorlinks]{hyperref}
\usepackage{latexsym}
\usepackage{amsfonts}
\usepackage{amssymb}
\usepackage{amsthm}
\usepackage{paralist}
\usepackage{mathrsfs}
\begin{document}
\newtheorem{thm}{Theorem}
\newtheorem{cor}[thm]{Corollary}
\newtheorem{lem}{Lemma}
\theoremstyle{remark}\newtheorem{rem}{Remark}
\theoremstyle{definition}\newtheorem{defn}{Definition}

\title{Oscillation inequalities on real and ergodic $H^1$ spaces}
\author{Sakin Demir}
\author{Sakin Demir}
\author{Sakin Demir\\
Agri Ibrahim Cecen University\\ 
Faculty of Education\\
Department of Basic Education\\
04100 A\u{g}r{\i}, Turkey\\
E-mail: sakin.demir@gmail.com
}

\maketitle

\renewcommand{\thefootnote}{}



\renewcommand{\thefootnote}{\arabic{footnote}}
\setcounter{footnote}{0}

\begin{abstract}Let $(x_n)$ be a sequence and $\rho\geq 1$. For a fixed sequences $n_1<n_2<n_3<\dots$, and $M$ define the oscillation operator
$$\mathcal{O}_\rho (x_n)=\left(\sum_{k=1}^\infty\sup_{\substack{n_k\leq m< n_{k+1}\\m\in M}}\left|x_m-x_{n_k}\right|^\rho\right)^{1/\rho}.$$
Let $(X,\mathscr{B} ,\mu , \tau)$ be a dynamical system with  $(X,\mathscr{B} ,\mu )$ a probability space and $\tau$ a measurable, invertible, measure preserving point transformation from $X$ to itself.\\
Suppose that the sequences $(n_k)$ and $M$ are lacunary. Then we prove the following results for $\rho\geq 2$:
\begin{enumerate} [(i)]
\item   Define $\phi_n(x)=\frac{1}{n}\chi_{[0,n]}(x)$ on $\mathbb{R}$. Then there exists a constant $C>0$ such that
$$\|\mathcal{O}_\rho (\phi_n\ast f)\|_{L^1(\mathbb{R})}\leq C\|f\|_{H^1(\mathbb{R})}$$
\noindent
 for all $f\in H^1(\mathbb{R})$.
\item  Let
$$A_nf(x)=\frac{1}{n}\sum_{k=1}^nf(\tau^kx)$$
be the usual ergodic averages in ergodic theory. Then 
 $$\|\mathcal{O}_\rho (A_nf)\|_{L^1(X)}\leq C\|f\|_{H^1(X)}$$
\noindent
for all $f\in H^1(X)$.
\item If $[f(x)\log (x)]^+$ is integrable, then $\mathcal{O}_\rho (A_nf)$ is integrable.
\end{enumerate}
\end{abstract}

\section{Preliminaries}
Let $(X,\mathscr{B} ,\mu )$ a totally $\sigma$-finite measure space and $\tau :X\to X$ be an ergodic measure preserving transformation. The function
$$f^{\ast}(x)=\sup_{n}\frac{1}{n}\sum_{i=0}^{n-1}\left |f(\tau^i x)\right |$$
is known as ergodic maximal function analogue to the  Hardy-Littlewood maximal function $Mf$ on the real line $\mathbb{R}$ given by 
$$Mf(x)=\sup_{x\in I}\frac{1}{|I|}\int_{I}\left|f(t)\right|\,dt,$$
where $I$ denotes an arbitrary interval in $\mathbb{R}$. \\

Let now $f$ be an integrable function and define
$$f^{\sharp}(x)=\sup_{n}\frac{1}{n}\sum_{i=0}^{n-1}\left|f(\tau^ix)-T_{n}f(x)\right|$$
where
$$T_nf(x)=\frac{1}{n}\sum_{i=0}^{n-1}\left|f(\tau^ix)\right|.$$
Now recall that the space $H^1$ on the real line $\mathbb{R}$ can  be characterized by
$$H^1(\mathbb{R})=\left\{f\in L^1(\mathbb{R}):\widetilde{H}f\in L^1(\mathbb{R})\right\}$$
with the norm 
$$\|f\|_{H^1}\sim\|f\|_1+\|\widetilde{H}f\|_1.$$
where $\widetilde{H}f$ is the Hilbert transform on $\mathbb{R}$ defined by
$$\widetilde{H}f(x)=\lim_{\epsilon\to 0^+}\int_\epsilon^\infty\frac{f(t+t)-f(x-t)}{t}\,dt.$$
Similar to the characterization by the maximal function we also define ergodic $H^1$ space by
$$H^1(X)=\left\{f\in L^1(X):Hf\in L^1(X)\right\}$$
with the norm 
$$\|f\|_{H^1}\sim\|f\|_1+\|Hf\|_1.$$
where $Hf$ is the ergodic Hilbert transform defined by
$$Hf(x)=\sum_{k=1}^{\infty}\frac{f(\tau^kx)-f(\tau^{-k}x)}{k}.$$
Similar to the classical case we can also identify the dual of ergodic $H^1$ spaces as ergodic bounded mean oscillation EBMO defined by the space of functions $f$ for which $f^{\sharp}$ is bounded with  EBMO norm given by 
$$\|f\|_{\rm{EBMO}}=\left\|f^{\sharp}\right\|_{\infty}$$
Let now $B$ be a Banach space and $p<\infty$, and let $f$ be a  $B$-valued (strongly) measurable function defined on $\mathbb{R}$. Then the Bochner-Lebesgue space $L_B^p=L_B^p(\mathbb{R})$ is defined as 
$$L_B^p=\left\{f:\|f\|_{L_B^p}<\infty\right\}$$
where
$$\|f\|_{L_B^p}=\left(\int_{\mathbb{R}}\|f(x)\|_B^p\,dx\right)^{1/p}.$$
When $B$ is the scalar field, we simply write $L^p$ and $\|\cdot\|_p$. We also define the space $WL_B^p={\rm{weak}}-L_B^p$ as the space of all $B$-valued functions $f$ such that
$$\|f\|_{WL_B^p}=\sup_{\lambda >0}\lambda\left(m\left\{x\in\mathbb{R}:\|f(x)\|_{B}\right\}>\lambda \right)^{1/p}<\infty.$$
When we replace Lebesgue measure by $w(x)dx$ for some positive weight $w$ in $\mathbb{R}$ we denote the corresponding spaces by $L_B^p(w)$ and $WL_B^p(w)$. When $p=\infty$, we write
$$L^{\infty}(B)=\left\{f:\|f\|_{L^{\infty}(B)}<\infty\right\},$$
where
$$\|f\|_{L^{\infty}(B)}=\textrm{ess}\sup\|f\|_{B}$$
and the space of all compactly supported members of $L^{\infty}(B)$ will  be denoted by $L_c^{\infty}(B)$.\\
For a locally integrable  $B$-valued function $f$, we define the maximal functions
$$M_rf(x)=\sup_{x\in I}\left(\frac{1}{|I|}\int_{I}\|f(y)\|_{B}^{r}\,dy\right)^{1/r},\quad 1\leq r\leq\infty ,$$
and
$$f^{\sharp}(x)=\sup_{x\in I}\frac{1}{|I|}\int_{I}\|f(y)-f_I\|_{B}\,dy,$$
where $I$ denotes an arbitrary interval in $\mathbb{R}$ and
$$f_I=\frac{1}{|I|}\int_If(t)\,dt$$
which is an element of $B$.\\
Note that $f^{\sharp}$ is the sharp maximal function in the classical case when $B=\mathbb{R}$ and $\|\cdot\|_B=|\cdot |$, $M_1f$ is the Hardy-Littlewood maximal function and $M_{\infty}f$ is the constant function.
Similar to the classical case we define the  $B$-valued BMO space as
$${\rm{BMO}}(B)=\left\{f\in L^1_{{\rm{loc.}}B}:\|f\|_{{\rm{BMO}}(B)}=\|f^{\sharp}\|_{L^\infty (B)}<\infty\right\}.$$
Given a  $B$-valued function $f$, we obtain a non-negative function $\|f\|_B$ defined by
$$\|f\|_B(x)=\|f(x)\|_B,$$
and it is important to point out that 
$$\left\|\left(\|f\|_B\right)\right\|_{{\rm{BMO}}}\leq 2\|f\|_{{\rm{BMO}}(B)}.$$
As usual a  $B$-atom is a function $a\in L^{\infty}(B)$ supported in an interval $I$ and such that
$$\|a(x)\|_B\leq\frac{1}{|I|},\quad \int_Ia(x)\,dx=0$$
and the space $H_B^1(\mathbb{R})$ such that
$$f(x)=\sum_j\lambda_ja_i(x);\quad (\lambda_j)\in l^1,$$
where $a_j$ are $B$-atoms with 
$$\|f\|_{H_B^1}=\inf\sum_j|\lambda_j|.$$
Similar to the classical case given $B\in\textrm{UMD}$ we also have
$$H_B^1(\mathbb{R})=\big\{f\in L_B^1(\mathbb{R}):\widetilde{H}f\in L_B^1(\mathbb{R})\big\}$$
and
$$\|f\|_{H^1_B}\sim\|f\|_{L_B^1}+\big\|\widetilde{H}f\big\|_{L_B^1}.$$


\begin{defn}\label{lacunary} A sequence $(n_k)$ of integers is called lacunary if there is a constant $\alpha >1$ such that
$$\frac{n_{k+1}}{n_k}\geq\alpha$$
for all $k=1,2,3,\dots$ .
\end{defn}

\begin{defn}Let $(x_n)$ be a sequence and $\rho\geq 1$. For a fixed sequences $n_1<n_2<n_3<\dots$, and $M$ define the oscillation operators
$$\mathcal{O}_\rho (x_n)=\left(\sum_{k=1}^\infty\sup_{\substack{n_k\leq m< n_{k+1}\\m\in M}}\left|x_m-x_{n_k}\right|^\rho\right)^{1/\rho}$$
and
$$\mathcal{O}_\rho^{\prime} (x_n)=\left(\sum_{k=1}^\infty\sup_{\substack{n_k\leq n\leq m< n_{k+1}\\n, m\in M}}\left|x_m-x_n\right|^\rho\right)^{1/\rho}.$$
\begin{rem}\label{o2l2rem}It is obvious that
$$\mathcal{O}_2(x_n)\leq \mathcal{O}_2^\prime (x_n)\leq 2\mathcal{O}_2(x_n).$$
\end{rem}
\end{defn}
Let $(X,\mathscr{B} ,\mu ,\tau)$ be a dynamical system with  $(X,\mathscr{B} ,\mu )$ a probability space and $\tau$ a measurable, invertible, measure preserving point transformation from $X$ to itself. Let $\rho\geq 2$ and $f\in L^1(X)$. Recall the usual ergodic averages
$$A_nf(x)=\frac{1}{n}\sum_{k=1}^nf(\tau^kx).$$
The we can define the oscillation operator
$$\mathcal{O}_\rho (A_nf)(x)=\left(\sum_{k=1}^\infty\sup_{\substack{n_k\leq m< n_{k+1}\\m\in M}}\left|A_mf(x)-A_{n_k}f(x)\right|^\rho\right)^{1/\rho}.$$
\begin{rem}Let $\rho\geq 1$. For a fixed sequences $n_1<n_2<n_3\dots$, and $M$, we can construct a Banach space $\mathbb{B}$ as an $\ell^\rho$ sum of finite-dimensional $\ell^\infty$ spaces with the following norm:\\
$$\|b\|_{\mathbb{B}}=\left(\sum_{k=1}^\infty\sup_{\substack{n_k\leq m< n_{k+1}\\m\in M}}\left|b_m-b_{n_k}\right|^\rho\right)^{1/\rho}.$$
\end{rem}
\noindent
Let us now define the kernel operator $K:\mathbb{R}\to \mathbb{B}$ as
$$K(x)=\left(\left(\frac{1}{m}\chi_{[0,m]}(x)-\frac{1}{n_k}\chi_{[0,n_k]}(x):n_k\leq m< n_{k+1},\: m\in M\right):k\geq 1\right).$$\\
Let 
 $\phi_n(x)=\frac{1}{n}\chi_{[0,n]}(x)$ and define 
$$\mathcal{O}_\rho (\phi_n\ast f)(x)=\left(\sum_{k=1}^\infty\sup_{\substack{n_k\leq m< n_{k+1}\\m\in M}}\left|\phi_m\ast f(x)-\phi_{n_k}\ast f(x)\right|^\rho\right)^{1/\rho}.$$
It is clear that
$$\mathcal{O}_\rho (\phi_n\ast f)(x)=\|K\ast f(x)\|_{\mathbb{B}},$$
where
$$K\ast f(x)=\int K(x-y)\cdot f(y)\, dy.$$
\section{The Results}
Our first result is the following theorem which plays a key role in proving inequalities for the oscillation operators $\mathcal{O}_\rho (\phi_n\ast f)$ and $\mathcal{O}_\rho (A_nf)$:
\begin{thm}\label{hconk}Let $(n_k)$ and $M$ be lacunary sequences, then there exists a constant $C>0$ such that
$$\int_{|x|>4|y|}\|K(x-y)-K(x)\|_{\mathbb{B}}\, dx\leq C$$
where $C$ does not depend on $y\in\mathbb{R}$, i.e., $K$ satisfies the Hörmander condition.
\end{thm}
\begin{proof}
\indent
For an $n\in\mathbb{N}$ let 
$$\phi_n(x)=\frac{1}{n}\chi_{[0,n]}(x)$$
and
$$\Phi_{(m,k)}(x)=\phi_m(x)-\phi_{n_k}(x).$$
Then we have
$$\|K(x-y)-K(x)\|_{\mathbb{B}}=\left(\sum_{k=1}^\infty\sup_{\substack{n_k\leq m< n_{k+1}\\m\in M}}\left|\Phi_{(m,k)}(x-y)-\Phi_{(m,k)}(x)\right|^\rho\right)^{1/\rho}.$$
Note that maximum is taken over all $m\in M$ between the related intervals as described in the definition of $\mathcal{O}_\rho (\phi_n\ast f)$.\\
Let us first consider the case $x>4y$, $y>0$. We have
\begin{align*}
\Phi_{(m,k)}(x-y)-\Phi_{(m,k)}(x)&=\phi_m(x-y)-\phi_{n_k}(x-y)-(\phi_m(x)-\phi_{n_k}(x))\\
                                                     &=\phi_m(x-y)-\phi_m(x)-(\phi_{n_k}(x-y)-\phi_{n_k}(x)).
\end{align*}
Since $x>4y$, $y>0$ we get
\begin{align*}
\phi_m(x-y)-\phi_m(x)&=\frac{1}{m}\chi_{[0,m]}(x-y)-\frac{1}{m}\chi_{[0,m]}(x)\\
                                   &=\frac{1}{m}\chi_{[y,m+y]}(x)-\frac{1}{m}\chi_{[0,m]}(x)\\
                                   &=\frac{1}{m}\chi_{[4y,m+y]}(x)-\frac{1}{m}\chi_{[4y,m]}(x)\\
                                   &=\frac{1}{m}\chi_{[m,m+y]}(x).
\end{align*}
Similarly, we have
$$\phi_{n_k}(x-y)-\phi_{n_k}(x)=\frac{1}{n_k}\chi_{[n_k,n_k+y]}(x).$$
We have
\begin{align*}
\int_{x>4y}\|K(x-y)-K(x)\|_{\mathbb{B}}\, dx= \;\;\;\;\;\;\;\;\;\;\;\;\;\;\;\;\;\;\;\;\;\;\;\;\;\;\;\;\;\;\;\;\;\;\;\;\;\;\;\;\;\;\;\;\;\;\;\;\;\;\;\;\;\;\;\;\;\;\;\\
=\int_{x>4y}\left(\sum_{k=1}^\infty\sup_{\substack{n_k\leq m< n_{k+1}\\m\in M}}\left|\frac{1}{m}\chi_{[m,m+y]}(x)-\frac{1}{n_k}\chi_{[n_k,n_k+y]}(x)\right|^\rho\right)^{1/\rho}\, dx\\
\leq  \int_{x>4y}\left(\sum_{k=1}^\infty\sup_{\substack{n_k\leq m< n_{k+1}\\m\in M}}\left|\frac{1}{m}\chi_{[m,m+y]}(x)\right|^\rho\right)^{1/\rho}\, dx+\\
+\int_{x>4y}\left(\sum_{k=1}^\infty\sup_{\substack{n_k\leq m< n_{k+1}\\m\in M}}\left|\frac{1}{n_k}\chi_{[n_k,n_k+y]}(x)\right|^\rho\right)^{1/\rho}\, dx\\
\leq  \int_{x>4y}\sum_{y\leq m}\sup_{\substack{n_k\leq m< n_{k+1}\\m\in M}}\frac{1}{m}\chi_{[m,m+y]}(x)\, dx+\\
+\int_{x>4y}\sum_{y\leq n_k}\sup_{\substack{n_k\leq m< n_{k+1}\\m\in M}}\frac{1}{n_k}\chi_{[n_k,n_k+y]}(x)\, dx\\
\leq   y\sum_{y\leq m}\frac{1}{m}+ y\sum_{y\leq n_k}\frac{1}{n_k}.\\
\end{align*}
On the other hand, we know that $(n_k)$ is a lacunary sequence, there is a constant $\beta$ such that
$$\frac{n_{k+1}}{n_k}\geq\beta >1.$$
Therefore, there is a constant $C(\beta )$ such that
$$\sum_{y\leq n_k}\frac{1}{n_k}\leq\frac{C(\beta ) }{y}.$$
Similarly, since $M$ is a lacunary sequence there is a constant $\alpha >1$  and a constant $C(\alpha )$ such that
$$\sum_{y\leq n_k}\frac{1}{m}\leq\frac{C(\alpha ) }{y}.$$
As a conclusion we see that
$$\int_{x>4y}\|K(x-y)-K(x)\|_{\mathbb{B}}\, dx\leq C(\alpha )+C(\beta )=C(\alpha ,\beta )$$
and this shows that the kernel operator $K$ satisfies the H\"ormander condition.\\
\indent
Let us now consider the case $y\leq 0$ and $x>4|y|$, 
$$\phi_m(x-y)-\phi_m(x)=\frac{1}{m}\chi_{[y,y+m]}(x)-\frac{1}{m}\chi_{[0,m]}(x)$$
we then have
$$\phi_m(x-y)-\phi_m(x)=\frac{1}{m}\chi_{[m+y,m]}(x)$$
and similarly, we have
$$\phi_{n_k}(x-y)-\phi_{n_k}(x)=\frac{1}{n_k}\chi_{[n_k+y, n_k]}(x)$$
for $x>4|y|$.\\
Again we see as in the previous case that
$$\int_{x>4y}\|K(x-y)-K(x)\|_{\mathbb{B}}\, dx\leq C(\alpha )+C(\beta )=C(\alpha ,\beta )$$
and this shows that $K$ satisfies the the H\"ormander condition in this case as well.\\
Suppose now that $x<0$ and $y>0$. Since $|x|>4|y|$, we see that $x-y<0$ thus in this case
$$\phi_m(x-y)-\phi_m(x)=\frac{1}{m}\chi_{[0,m]}(x-y)-\frac{1}{m}\chi_{[0,m]}(x)=0$$
and
$$\phi_{n_k}(x-y)-\phi_{n_k}(x)=\frac{1}{n_k}\chi_{[0,n_k]}(x-y)-\frac{1}{n_k}\chi_{[0,n_k]}(x)=0.$$
Thus for any $y$ we have 
$$\int_{|x|>4|y|}\|K(x-y)-K(x)\|_{\mathbb{B}}\, dx=0.$$
We finally need to consider the case $x<0$ and $y<0$. In this case we have
$$\phi_m(x-y)-\phi_m(x)=\frac{1}{m}\chi_{[y, y+m]}(x)$$
and
$$\phi_{n_k}(x-y)-\phi_{n_k}(x)=\frac{1}{n_k}\chi_{[y, y+n_k]}(x).$$
Since we also have $|x|>4|y|$, we see that
$$\phi_m(x-y)-\phi_m(x)=\frac{1}{m}\chi_{[y, y+m]}(x)=0$$
and 
$$\phi_{n_k}(x-y)-\phi_{n_k}(x)=\frac{1}{n_k}\chi_{[y, y+n_k]}(x)=0.$$
Thus for any $y$ we have 
$$\int_{|x|>4|y|}\|K(x-y)-K(x)\|_{\mathbb{B}}\, dx=0.$$
\end{proof}
\begin{lem}\label{ol2lem}Suppose that $(n_k)$ is a lacunary sequence. Then there exits a constant $C>0$ such that
$$\|\mathcal{O}_2(\phi_n\ast f)\|_{ L^2(\mathbb{R})}\leq C\|f\|_{ L^2(\mathbb{R})}$$
for all $f\in L^2(\mathbb{R})$.
\end{lem}
\begin{proof}When $(n_k)$ is lacunary it is proven in R. L. Jones {\it{et al}} \cite{rjkm}) (see Corollary 2.11) that there exists a constant $C>0$ such that
$$\|\mathcal{O}_2^\prime (A_nf)\|_{ L^2(X)}\leq C\|f\|_{ L^2(X)}$$
for all $f\in L^2(X)$.  Remark~\ref{o2l2rem} clearly implies that 
$$\|\mathcal{O}_2(A_nf)\|_{ L^2(X)}\leq \|\mathcal{O}_2^\prime (A_nf)\|_{ L^2(X)}.$$
Thus we get
$$\|\mathcal{O}_2(A_nf)\|_{ L^2(X)}\leq C\|f\|_{ L^2(X)}$$
$f\in L^2(X)$. The proof of our lemma follows by an application of Calder\'on transfer principle.
\end{proof}
\begin{thm}\label{l2of} There exists a constant $C>0$ such
$$\|\mathcal{O}_\rho (\phi\ast f)\|_{ L^2(\mathbb{R})}\leq C\|f\|_{ L^2(\mathbb{R})}$$
for all $f\in L^2(\mathbb{R})$.
\end{thm}
\begin{proof} Let $\rho\geq2$. Then we obtain
\begin{align*}
\|\mathcal{O}_\rho (\phi\ast f)\|_{ L^2(\mathbb{R})}^2&=\int\left(\sum_{k=1}^\infty\sup_{\substack{n_k\leq m< n_{k+1}\\m\in M}}\left|\phi_n\ast f(x)-\phi_{n_k}\ast f(x)\right|^\rho\right)^{2/\rho}\, dx\\
&\leq \int\sum_{k=1}^\infty\sup_{\substack{n_k\leq m< n_{k+1}\\m\in M}}\left|\phi_n\ast f(x)-\phi_{n_k}\ast f(x)\right|^2\, dx\\
&=\int\left[\left(\sum_{k=1}^\infty\sup_{\substack{n_k\leq m< n_{k+1}\\m\in M}}\left|\phi_n\ast f(x)-\phi_{n_k}\ast f(x)\right|^2\right)^{1/2}\right]^2\, dx\\
&=\|\mathcal{O}_2 (\phi\ast f)\|_{ L^2(\mathbb{R})}^2\\
&\leq C\|f\|_{ L^2(\mathbb{R})}^2
\end{align*}
for some constant $C>0$ by Lemma~\ref{ol2lem}.
\end{proof}
We can now state and prove our next result:
\begin{thm}\label{L1ToH1R}There exits a constant $C>0$ such that
$$\|\mathcal{O}_\rho(\phi_n\ast f)\|_{L^1(\mathbb{R})}\leq C\|f\|_{H^1(\mathbb{R})}$$
for all $f\in H^1(\mathbb{R})$.
\end{thm}
\begin{proof}It suffices to show that
$$\|\mathcal{O}_{\rho}(\phi_n\ast a)\|_{L^1(\mathbb{R})}\leq C$$
for any atom $a$ with constant $C$ independent of the choice of $a$. We first consider a $1$-atom centered at $0$ with support $\text{support}(a)\subset I_R$, where $I_R$ denotes an interval centered at $0$ with side length $2R$.
Since 
\begin{align*}
\|a\|^2_{L^2(\mathbb{R})}&=\int_{\mathbb{R}}|a(x)|^2\, dx\\
&=\int_{I_R}|a(x)|^2\, dx\\
&\leq \frac{1}{|I_R|^2}\int_I\, dx\\
&=\frac{1}{|I_R|}.
\end{align*}
we see that
$$\|a\|_{L^2(\mathbb{R})}\leq \frac{1}{|I_R|^{1/2}}.$$
Because of the definition of an atom we also have
$$\int_{I_R}a(x)\, dx=0.$$
Hence
\begin{align*}
\int_{|x|\geq 4R}|\mathcal{O}_\rho (\phi_n\ast a)(x)|\, dx&=\int_{|x|\geq 4R}\|K\ast a(x)\|_{\mathbb{B}}\, dx\\
&=\int_{|x|\geq 4R}\left|\int_{I_R}\|\{K(x-y)-K(x)\}\cdot a(y)\|_{\mathbb{B}}\, dy\right|\, dx\\
&\leq \int_{I_R}|a(y)|\, dy\int_{|x|\geq 4|y|}\|K(x-y)-K(x)\|_{\mathbb{B}}\, dx\\
&\leq C(\alpha , \beta )\int_{I_R}|a(y)|\, dy\\
&\leq C(\alpha , \beta ).
\end{align*}
On the other hand by Hölder's inequality
\begin{align*}
\int_{|x|< 4R}|\mathcal{O}_\rho (\phi_n\ast a)(x)|\, dx&=\int_{|x|< 4R}\|K\ast a(x)\|_{\mathbb{B}}\, dx\\
&\leq \left(\int \|K\ast a(x)\|_{\mathbb{B}}^2\right)^{1/2}2|I_R|^{1/2}\\
&\leq C_1\|a\|_{L^2(\mathbb{R})}2|I_R|^{1/2}\\
&\leq C.
\end{align*}
We obtain
$$\|\mathcal{O}_{\rho}(\phi_n\ast a)\|_{L^1(\mathbb{R})}\leq C$$
in both inequalities for any atom $a$ centered at the origin.\\
Let now $b$ be an atom centered at $c\in \mathbb{R}$. Then $a(x)=b(x-c)$ is an atom centered at $0$. Moreover, we have
$$\|\mathcal{O}_{\rho}(\phi_n\ast a)\|_{L^2(\mathbb{R})}\leq C_1 \|a\|_{L^2(\mathbb{R})}$$
and we  obtain as before
$$\|\mathcal{O}_{\rho}(\phi_n\ast a)\|_{L^1(\mathbb{R})}\leq C.$$
Hence we have
\begin{align*}
\|\mathcal{O}_{\rho}(\phi_n\ast b)\|_{L^1(\mathbb{R})}&=\|\mathcal{O}_{\rho}(\phi_n\ast a)\|_{L^1(\mathbb{R})}\\
&\leq C.
\end{align*}
\end{proof}

When we consider the maximal function characterization of the space $H^1(\mathbb{R})$, it is clear that $\|f\|_{H^1(\mathbb{R})}$ is the $L^1$ norm of an operator of convolution type. That's because we can apply the transfer principle of S. Demir~\cite{demir} to Theorem~\ref{L1ToH1R} to find the following result with the same constant as in Theorem~\ref{L1ToH1R}:
\begin{thm}\label{L1ToH1X}There exits a constant $C>0$ such that
$$\|\mathcal{O}_\rho (A_nf)\|_{L^1(X)}\leq C\|f\|_{H^1(X)}$$
for all $f\in H^1(X)$.
\end{thm}
Note that in Theorem~\ref{L1ToH1X} $H^1(X)$ denotes the ergodic Hardy space. We refer those readers who are not so familiar with these spaces to  R.~Caballero and A.~de la Torre~\cite{cdelatorre}.\\
Recall the following well known theorem of D.~Ornstein~\cite{dornstn} which will be used in the proof of our last result:
\begin{lem}\label{lloglmax}The ergodic maximal function
$$f^{\ast}(x)=\sup_{n}\frac{1}{n}\sum_{i=0}^{n-1}\left |f(\tau^i x)\right |$$
is integrable if and only if  $[f(x)\log (x)]^+$ is integrable, where $g^+$ denotes the positive part of $g$.
\end{lem}
Our last result gives a condition for integrability of the oscillation operator $\mathcal{O}_\rho (A_nf)$.
\begin{thm}\label{intos}If $[f(x)\log (x)]^+$ is integrable, then $\mathcal{O}_\rho (A_nf)$ is integrable for $\rho\geq 2$.
\end{thm}
\begin{proof}It is proven in S. Demir~\cite{sdem} that
$$\|f\|_{H^1(X)}\leq \|f^\ast\|_{L^1(X)}.$$
Now suppose that $[f(x)\log (x)]^+$ is integrable, then Lemma~\ref{lloglmax} implies the finiteness of $\|f^\ast\|_{L^1(X)}$, and thus the finiteness of $\|f\|_{H^1(X)}$. This means that $f\in H^1(X)$, and by Theorem~\ref{L1ToH1X} there exists a constant $C>0$ such that
$$\|\mathcal{O}_\rho (A_nf)\|_{L^1(X)}\leq C\|f\|_{H^1(X)},\;\;\;\rho \geq 2$$
thus we have
$$\|\mathcal{O}_\rho (A_nf)\|_{L^1(X)}<\infty$$
and this completes our proof.
\end{proof}


\begin{thebibliography}{99}
\bibitem{cdelatorre}R.~Caballero and A.~de la Torre, 
\emph{An atomic theory of ergodic $H^p$ spaces},
Studia Math. 82 (1985) 39-69.
\bibitem{apcal}A. P.~Calder\'on, 
\emph{Ergodic theory and translation-invariant operators},
Proc. Nat. Acad. Sci. USA. 59 (1968) 349-353.
\bibitem{demir}S.~Demir,
\emph{A generalizaition of Calder\'on transfer principle},
Journal of Computer \& Mathematical Sciences 9(5) (2018) 325-329.
\bibitem{sdem}S.~Demir, 
\emph{$H^p$ spaces and inequalities in ergodic theory},
Ph.D Thesis,  University of Illinois at Urbana-Champaign, USA, May 1999.
\bibitem{rjkm} R. L.~Jones, R.~Kaufman, J. M.~Rosenblatt and Máté Wierdl,
\emph{Oscillation in ergodic theory},
Ergodic Th. \& Dynam. Sys 18 (1998) 889-935.                                                                                                                                                                                                                                                                        
\bibitem{dornstn}D.~Ornstein, 
\emph{A remark on Birkoff ergodic theorem},
Illinois J. Math. 15 (1971) 77-79.
\end{thebibliography}
\end{document}